\documentstyle[12pt, righttag, amssym]{amsart}
\def \Z{\Bbb Z}
\def \M{\Bbb M}
\def \C{\Bbb C}
\def \R{\Bbb R}
\def \Q{\Bbb Q}

\def \wt{{\rm wt}}
\def \tr{{\rm tr}}

\def \End{{\rm End}}

\def \Aut{{\rm Aut}}

\def \mod{{\rm mod}}

\def \<{\langle}
\def \>{\rangle}

\def \a{\alpha }

\def \L{\Lambda }

\def \o{\omega }

\def \V{V^{\natural}}
\def \voa{vertex operator algebra\ }
\def \be{\begin{equation}}
\def \ee{\end{equation}}
\def \bl{\begin{lem}}
\def \el{\end{lem}}
\def \ba{\begin{array}}
\def \ea{\end{array}}
\def \bt{\begin{thm}}
\def \et{\end{thm}}
\def \ch{{\rm ch}}
\newtheorem{th1}{Theorem}
\newtheorem{ree}[th1]{Remark}
\newtheorem{thm}{Theorem}[section]

\newtheorem{coro}[thm]{Corollary}
\newtheorem{lem}[thm]{Lemma}

\numberwithin{equation}{section}

\begin{document}
\title[Some twisted sectors for the Moonshine Module]{Some twisted sectors for
the Moonshine Module}
\author{Chongying Dong, Haisheng Li and Geoffrey Mason}
\address{Department of Mathematics, University of California, Santa Cruz}
    \email{dong@@cats.ucsc.edu\\
    hli@@cats.ucsc.edu\\
    gem@@cats.ucsc.edu}
    \thanks{This paper is in final form and no version of it will
be submitted for publication elsewhere.}
    \thanks{C.D. and G.M. are partially supported
by NSF grants and a research grant from the Committee on Research, US Santa
Cruz.}
    \subjclass{Primary 17B69; Secondary 17B68, 81T40}
	\bibliographystyle{alpha}
	\maketitle
\begin{abstract}
The construction of twisted sectors, or $g$-twisted modules, for a
\voa  $V$ and automorphism $g,$ is a fundamental problem in algebraic
conformal field theory and the theory of orbifold models. For the moonshine
module $V^{\natural},$ whose automorphism is the Monster ${\Bbb M},$ Tuite
has shown that this problem is intimately related to the generalized moonshine
conjecture which relates hauptmoduln to twisted sectors.

In this paper we show how to give a uniform existence proof for
irreducible $g$-twisted $V^{\natural}$-modules for elements of type
$2A,$ $2B$ and $4A$ in $\M.$ The most interesting of these is the
twisted sector $\V(2A),$ whose automorphism group is essentially
the centralizer of $2A$ in $\M.$ This is a 2-fold central extension
of the Baby Monster, the second largest sporadic simple group.

We also establish uniqueness of the twisted sectors and the hauptmodul
property for the graded traces of automorphisms of odd order, as
predicted by conformal field theory.
\end{abstract}

\section{Introduction}
\setcounter{equation}{0}
The purpose of this paper is to give a uniform existence proof for irreducible
$g$-twisted $\V$-modules for elements $g$ in $\M$ of type $2A,2B$ and $4A.$
Here $\V$ is the moonshine module [FLM] with automorphism group the Monster
$\M,$ and the notation for elements in $\M$ follows the conventions of [Cal].
We assume the reader to be familiar with the definition of {\em vertex operator
algebra} (VOA) (see [B1] and [FLM]) which we do not repeat here, noting only
that $V^{\natural}$ is a particularly famous example of a VOA.

If $V$ is a VOA and $g$ an automorphism of finite order, there is the notion
of a ({\em weak}) $g$-{\em twisted} $V$-{\em module}. Briefly, a weak
$g$-twisted $V$-module is a pair $(M,Y_M)$ consisting of a $\C$-graded linear
space
$M=\oplus_{n\in\C}M_{n}$ locally truncated below in the sense that
$M_{n+q}=0$ for fixed $n\in \C$ and sufficiently small $q\in \Q,$
 together with a linear map
$$\begin{array}{l}
V \to (\mbox{End}\,M)[[z^{1/T},z^{-1/T}]]\\
v\mapsto Y_M(v,z)=\displaystyle{\sum_{n\in\frac{1}{T}{\Z}}v_nz^{-n-1}}\ \ \ \
(v_n\in
\mbox{End}\,M),
\end{array}$$
where $T$ is the order of $g.$ One requires a certain Jacobi identity
to be satisfied by these maps, as well as an appropriate action of the Virasoro
algebra. See [DM1] for the complete definition.
A $g$-{\em twisted} $V$-{\em module} is a weak $g$-twisted $V$-module such that
all of the homogeneous
spaces $M_n$ are of finite dimension. We should point out that the definition
in [DM1] differs slightly from that used here. Namely, we allow $M$ to be
$\C$-graded in the present paper (in order to be able to apply results in
[DM2]) and we admit the possibility of infinite-dimensional homogeneous
spaces.

It can be easily shown that any weak $g$-twisted module is a direct sum
of submodules of type
\begin{equation}\label{1.1}
N=\oplus_{n=0}^{\infty}N_{c+\frac{n}{T}}
\end{equation}
where $c\in \C$ such that $N_c\ne 0.$ The subspace $N_c$ is called the
{\em top level} of $N.$ Note that the grading
of a simple weak $g$-twisted module always has the form (\ref{1.1}).

\begin{th1}\label{t1} Let $g\in\M$ be of type $2A,2B$ or $4A.$ The following
hold:

(i) Every simple weak $g$-twisted $\V$-module is a $g$-twisted $\V$-module.

(ii) Up to isomorphism there is exactly one simple $g$-twisted $\V$-module.

(iii) If $g$ has type $2A$ then every $g$-twisted $\V$-module is completely
reducible.

\end{th1}

Twisted sectors for $\V$ of type $2B$ have been constructed by [Hu] and also
by two of us [DM2], but until now that has been the extent of our knowledge
of the existence of twisted sectors for $\V.$

Continuing with earlier notation, we define an {\em extended automorphism}
 of $(M,Y_M)$ to be a pair $(x,\a(x))$ where $x: M\to M$ and $\a(x):
V\to V$ are invertible linear maps satisfying
\begin{equation}\label{1.2}
\begin{array}{c}
xY_M(v,z)x^{-1}=Y_M(\alpha(x)v,z)\\
\alpha(x)g=g\alpha(x), \a(x){\bf 1}={\bf 1},
\a(x)\omega=\omega
\end{array}
\end{equation}
for $v\in V$ where ${\bf 1},\o$ are the vacuum and conformal element
respectively. This definition is a slight modification of [DM4], where
it is explained that if $V$ and $M$ are both simple then $x\mapsto
\alpha(x)$ is a group homomorphism from the group $\Aut^e(M)$ of
extended automorphisms of $M$ (which we identify with the group of
linear maps $x$) into $\Aut(V).$ Moreover, the kernel is a central
subgroup of $\Aut^e(M)$ consisting of the scalar operators.

We will show that if $\V(2A)$ is the simple $2A$-twisted $\V$-module whose
existence is guaranteed by Theorem \ref{t1}, then  $\Aut^e(\V(2A))$ is
isomorphic to $\C^*\times 2Baby$ where $2Baby$ denotes the 2-fold central
extension
of the Baby Monster which is known to be isomorphic to the centralizer of
$2A$ in $\M$ (cf. [Cal] and references therein). So in this case, the map
$x\mapsto \a(x)$ splits, though this is not always the case.
(E.g., for $g=2B$ it
does not. See Section 6 for details.) We identify $2Baby\subset \Aut^e(\V(2A))$
with the corresponding centralizer in $\M.$ If $V^{\natural}$ has grading
\be\label{1.3}
\V=\oplus_{n=0}^{\infty}\V_n
\end{equation}
then for $h\in \M$ we define
\begin{equation}\label{1.4}
Z(1,h,\tau)=q^{-1}\sum_{n=0}^{\infty}\tr(h|\V_n)q^n,
\end{equation}
sometimes called the {\em McKay-Thompson}
 series of $h.$ Here $q=e^{2\pi i\tau}$ as
usual. It is known (the Conway-Norton conjecture = Borcherds'
theorem [B2]) that each $Z(1,h,\tau)$ is a so-called hauptmodul.
For example, $Z(1,1,\tau)=J(q)=q^{-1}+196884q+\cdots$ is the absolute
modular invariant with constant term 0.

It transpires that if $M=\V(2A)$ then the constant $c$ in (\ref{1.1})
is equal to $1/2.$ Then for $g\in \M$ of type of $2A$ and
$h\in C_{\M}(g)\simeq 2Baby$ we define
\begin{equation}\label{1.5}
Z(g,h,\tau)=q^{-1}\sum_{n=1}^{\infty}\tr(h|\V(2A)_{n/2})q^{n/2}.
\end{equation}

Remark that from (\ref{1.2}) with $v$ equal to the conformal element $\o,$
it follows that $\Aut^e(\V(2A))$ preserves the $\Q$-grading on $\V(2A),$ so
that (\ref{1.5}) makes sense. We prove
\begin{th1}\label{t2} The following hold:

(i) If $h$ has odd order then
\begin{equation}\label{1.6}
Z(g,h,2\tau)=Z(1,gh,\tau).
\end{equation}

(ii) If $h$ either has odd order or satisfies $g\in\<h\>,$ then
$Z(g,h,\tau)$ is a hauptmodul.
\end{th1}

\begin{ree} 1. These results amount to establishing the generalized moonshine
conjecture due to Norton (cf. the appendix to [M]) for the commuting pairs
$(g,h)$ such that $\<g,h\>$ is cyclic and $g=2A.$

2. If $g\in\<h\>$ there are formulas which are analogous to (\ref{1.6}),
though more complicated.

3. Similar results also hold for the twisted sectors of type $2B$ and $4A.$

4. We can compute $Z(g,h,\tau)$ in other cases too. For example if
$\<g,h\>\simeq\Z_2\times \Z_2$ has all 3 involutions of type of $2A$
then we will show in Section 5
 that $Z(g,h,\tau)$ is precisely the hauptmodul $t_{2/2}=\sqrt{J(q)-984}$
as predicted by Conway-Norton [CN].
\end{ree}

The proofs of the two theorems depend heavily on the papers [DM2] and
[DMZ]. In [DMZ] it is explained that $\V$ contains a sub VOA which is
isomorphic to $L=L(\frac{1}{2},0)^{\otimes 48},$ the tensor product of
48 vertex operator algebras $L(\frac{1}{2},0)$ associated with the highest
weight unitary representation of the Virasoro algebra with central charge
$1/2.$ Moreover, $L$ is a rational VOA.

There are many sub VOAs of $\V$ isomorphic to $L;$ they may be
constructed from 48-dimensional associative subalgebras $A$ of the
Griess algebra (cf. [MN] and [Mi] for example) which themselves are
related to what we have called marked Golay codes [CM]. One can show
easily that if $g\in\M$ has type $2A, 2B$ or $4A$ then $g$ fixes a
suitable $A$ (element-wise), so that we may take $L$ to lie in the
sub VOA $(\V)^{\<g\>}$ of $g$-invariants.
Then we are in a position to invoke several results
from [DM2] to conclude Theorem \ref{t1}.

Theorem \ref{t1} is purely an existence theorem. To get further
information one needs to know the representation of the extended
automorphism group on the top level $M_c$ of the corresponding twisted
sector. This can vary greatly: for $M=\V(2B)$ the top level is of dimension
$24,$ whereas for $M=\V(2A)$ it is 1. In each case the extended
automorphism group
acts irreducibly. To establish Theorem \ref{t2} we argue first that
the top level of $\V(2A)$ is a {\em trivial}
 module for $2Baby,$ using the structure
of the Leech lattice as well as the fusion rules for $L(\frac{1}{2},0)$ and its
modules established in [DMZ]. Then we can prove that the dimension is 1 by
combining some monstrous calculation of Norton [N] together with
further properties of the Virasoro algebra.

At this point we are in a position to combine the modular-invariance
properties of [DM2], which relate graded traces of elements on $\V(2A)$ to
those on $\V,$ with Borcherds' proof of the original moonshine conjecture. This
leads to the proof of Theorem \ref{t2}.

We thank the referee for useful comments.

\section{Transpositions and Virasoro algebras}
\setcounter{equation}{0}

With the conformal grading (\ref{1.3}) of the moonshine module one knows that
$\V_0=\C {\bf 1}$ is spanned by the vacuum and $\V_1=0.$ $\V_2$ carries
the structure of a commutative non-associative algebra which we denote
by $B.$ As Monster module, $B$ is the direct sum of two irreducible
modules $\C\o\oplus B_0.$ We refer to both $B$ and $B_0$ as `the'
{\em Griess algebra.} Note that $\frac{1}{2}\o$ is the identity of $B.$

A {\em transposition} of $\M$ is an involution of type $2A.$ If $x$ is a
transposition then $C_{\M}(x)\simeq 2Baby.$ We fix such an $x$ and let
$H=C_{\M}(x).$ The space $B^H$ of $H$-invariants on $B$ is 2-dimensional and
spanned by $\o$ together with the so-called {\em transposition axis}
$t_x$ (cf. [C] or [MN]). Moreover $H$ is the subgroup of $\M$ leaving $t_x$
invariant. From this we easily conclude

\begin{lem}\label{l2.1}
Let $x_1,...,x_k$ be transpositions, let $A=\<t_{x_1},...,t_{x_k},\o\>$
be the subalgebra of $B$ spanned by $w$ and the corresponding
transpositions axes, and let $E=\<x_1,...,x_k\>$ be the subgroup of $\M$
generated by the $x_i.$ Then the subgroup $C_{\M}(E)$ of elements of $\M$
commuting with $E$ coincides with the subgroup $C_{\M}(A)$ of $\M$ fixing
$A$ (elementwise).\ $\Box$
\end{lem}

The subalgebra $A$ of Lemma \ref{l2.1} is associative if, and only if,
each product $x_ix_j$ ($i\ne j$) is an involution of type $2B$
(Theorem 5, Corollary 1 of [MN]). In particular, $E$ is elementary abelian
in this case.

It was shown in [MN] that we may choose an associative algebra
$A$ of dimension 48 satisfying the conclusions of Lemma \ref{l2.1}.
It is sufficient to give 24 mutually orthogonal vectors $v_i$ in the Leech
lattice $\Lambda$ such that $v_i$ has squared length 4, in which case the 48
elements $\pm v_i$ correspond to 48 mutually orthogonal transposition axes
of $B$ (cf. [DMZ] and [Mi]).

One has some latitude in choosing the $v_i.$ We may take, for example,
$v_1=(4,4,0^{22}),$ $v_2=(4,-4,0^{22}),$ $v_3=(0^2,4,4,0^{20}),$
$v_4=(0^2,4,-4,0^{20}),$. . . where the squared length of $v_i$ is
$\frac{1}{8}v_i\cdot v_i.$ Similarly, we could choose any twelve pairs of
coordinates and take the corresponding vectors. Such choices are not all
equivalent,
and in combination with a choice of Golay Code
give rise to the notion of a ``marked'' Golay Code [CM].

\begin{lem}\label{l2.2} We may choose the $v_i$ so that the corresponding
elementary abelian group $E$ satisfies $|E|\leq 2^9.$
\end{lem}

\pf We may take $E\leq Q=O_2(C)\simeq 2_+^{1+24}$ where $C$ is the
centralizer of an involution of type $2B$ in $M.$ We identify $Q/Z(Q)$
with $\L/2\L,$ so that each $\pm v_i$ becomes an involution (of type $2A$) in
$Q.$ Clearly $E=\<\pm v_i\>=\<-1,v_i\>.$ We have
$v_1+v_2=(8,0^{23})$ and for $i\geq 1,$ $v_{2i+1}+v_{2i+2}=(0^{2i},8,0^{23-2i})
=v_1+v_2+2(-4,0^{2i-1},4,0^{23-2i})\equiv v_1+v_2$ (\mod\,$2\L$). So
$E=\<-1,v_2,v_{2i-1},1\leq i\leq12\>.$

Suppose now that we mark our Golay Code so that the six 4-element sets
$\{1,2,3,4\},$ $...,$ $\{21,22,23, 24\}$
constitute a {\em sextet} i.e., the union
of any two of them is an {\em octad} (=block in the Witt design). Then
$v_1+v_3+v_5+v_7=(4^8,0^{16})=2(2^8,0^{16})\in2\L$ and similarly
$v_1+v_3+v_{4i+1}+v_{4i+3}\in2\L$ for $1\leq i\leq 5.$ Hence
$E=\<-1,v_2,v_1,v_3,v_5,v_9,
v_{13},$ $v_{17},v_{21}\>.$ \ \ $\Box$

\begin{coro}\label{c2.3} With the choice made in Lemma \ref{l2.2}, the
corresponding associative algebra $A$ of dimension 48 is fixed (pointwise) by
elements in $\M$ of types $2A,$ $2B$ and $4A.$
\end{coro}

\pf These are the three types of Monster elements (apart from 1) contained in
$Q.$ Since $Q$ is extra-special of order $2^{25}$ and $E\leq Q$ has order less
than or equal to $2^9,$ then certainly $C_Q(E)$ contains elements of these
three types. Now apply Lemma \ref{l2.1}.\ \ $\Box$

It was shown in [DMZ] that corresponding to the 48 transpositions
$x_1,...,x_{48}$ with transposition axes
spanning a 48-dimensional associative algebra we obtain
a certain Virasoro algebra of central charge 24. To explain this, first
recall that the Virasoro
algebra $Vir$ is spanned by $L_n,$ $n\in\Z$ and a central element $c$
satisfying  the relation
\begin{equation}\label{2.1}
[L_m,L_n]=(m-n)L_{m+n}+\frac{m^3-m}{12}\delta_{m+n,0}c.
\end{equation}

\begin{lem}\label{add} Let $\omega_1$ and $\omega_2$ be two orthogonal
idempotents such that the components of the
vertex operator $Y(\o_i,z)$ generate a copy Virasoro algebra
with central charge $1/2.$ Then the actions of these two Virasoro algebras
are commutative on $\V.$
\end{lem}

\pf Set $Y(\o_i,z)=\sum_{n\in\Z}L^i(n)z^{-n-2}$ for $i=1,2.$
 Then the inner product
of $\o_1$ with $\o_2$ is given by $L^1(2)\o_2,$ which is 0 by assumption.
By Norton's inequality (cf [C] or [MN])
$0=(\o_1,\o_2)=(\o_1^2,\o_2^2)\geq (\o_1\o_2,\o_1\o_2)
\geq 0,$ the product $\o_1\o_2$ (which is $L^1(0)\o_2$) in the
Griess algebra is also 0.
Thus $L^1(n)\o_2=0$ for all nonnegative integrals $n$
and $\o_2$ is a highest weight vector with highest weight 0 for the
Virasoro algebra $Vir_1$ generated by $L^1(m).$ The submodule of $\V$ for
$Vir_1$ generated by $\o_2$ is necessarily isomorphic to $L(1/2,0).$
Here $L(k,h)$ is the simple highest weight representation
of $Vir$ with central change $k$ and highest weight $h.$ From the module
structure of $L(1/2,0)$ we see immediately that $L^1(-1)\o_2=0.$ Now use
the commutator formula $[L^1(m),L^2(n)]=\sum_{i=-1}^{\infty}{m+1\choose i+1}
(L^1(i)\o_2)_{m+n+1-i}=0$ to complete the proof. \qed
\bigskip

Some multiple $\o_i$ of the transposition axis $t_{x_i}$ is an idempotent
such that the
components of the vertex operator $Y(\o_i,z)$ generate a copy Virasoro algebra
with central charge $1/2.$ As the $\o_i$ are orthogonal the algebras mutually
commute as operators on $\V$ by Lemma \ref{add}, yielding an action of the sub
VOA
\begin{equation}\label{2.2}
L=L(\frac{1}{2},0)^{\otimes 48}
\end{equation}
on $\V.$

Note that if $g\in\M$ fixes the associative algebra $A$ spanned by the
$t_{x_i}$ then $g$ fixes each $\o_i.$ Hence we get
\begin{lem}\label{l2.4}
If $g$ is of type $2A,$ $2B$ or $4A$ then the fixed subalgebra $(\V)^{\<g\>}$
contains $L.$
\end{lem}

The sum of all the
$\o_i$ is the conformal vector $\o.$ If $H=2Baby$ fixes $\o_1,$ say, then it
also fixes $\o_0$ where $\o=\o_1+\o_0.$ Now the component operators
of $Y(\o_0,z)$ generate a copy of the Virasoro algebra of central charge
$47/2,$ and in any case we get a subspace of $(\V)^H$ corresponding to
$Y(\o_1,z)$ and $Y(\o_0,z),$ that is to say the subspace
$L(\frac{1}{2},0)\otimes L(\frac{47}{2},0)$ which is again a vertex operator
algebra [FHL].

Now one knows (cf. [KR] for example) that the $q$-characters of these two
algebras are as follows:
\begin{equation}\label{2.3}
\begin{array}{c}
\displaystyle{
\ch_qL(\frac{1}{2},0)=\frac{1}{2}\left(\prod_{n=1}^{\infty}(1+q^{n-1/2})+
\prod_{n=1}^{\infty}(1-q^{n-1/2})\right)}\\[10pt]
\displaystyle{
\ch_qL(\frac{47}{2},0)=\prod_{n=2}^{\infty}(1-q^{n})^{-1}.}
\ea
\end{equation}
We thus calculate that the subspace $L(\frac{1}{2},0)\otimes L(\frac{47}{2},0)
\subset (\V)^H$ has $q$-character
\begin{equation}\label{2.4}
1+2q^2+2q^3+5q^4+6q^5+12q^6+\cdots.
\end{equation}

\begin{lem}\label{l2.5} Up to weight six, all elements of $(\V)^H$ lie in
$L(\frac{1}{2},0)\otimes L(\frac{47}{2},0).$
\end{lem}

\pf Norton has calculated (Table 5 of [N]) the decomposition of the permutation
representation of $\M$ on the conjugacy class of transpositions into
simple characters for $\M.$ The representation is multiplicity-free
and the characters occurring are precisely $\chi_1,\chi_2,\chi_4,$ $\chi_5,$
$\chi_9,
\cdots$ where we order the simple characters according to their degrees as in
[Cal].

By elementary character theory, if $V_i$ is the $\M$-module affording $\chi_i$
then the dimension of the space of $H$ (=$2Baby$) invariants $V_i^H$
is precisely the multiplicity of $\chi_i$ above. Hence it is either $0$ or
$1$ and is $1$ precisely if $i=1,2,4,5,9,\cdots.$

Now the decomposition of the first few homogeneous spaces $\V_n$ into
simple Monster module is known (e.g. [MS]). We have
$$\V_0=V_1,\ \ \V_1=0$$
$$\V_2=V_1\oplus V_2,\ \ \V_3=V_1\oplus V_2\oplus V_3$$
\begin{equation}\label{2.5}
\V_4=2V_1\oplus 2V_2\oplus V_3\oplus V_4
\end{equation}
$$\V_5=2V_1\oplus 3V_2\oplus 2V_3\oplus V_4\oplus V_6$$
$$\V_6=4V_1\oplus 5V_2\oplus 3V_3\oplus 2V_4\oplus V_5\oplus V_6\oplus V_7.$$

So to compute $\sum_{n=0}^6\dim(\V_n)^Hq^n,$ from the foregoing it is
equivalent to counting the number of occurrences of
$\chi_1,\chi_2,\chi_4,\chi_5$ in (\ref{2.5}). We find that the multiplicities
are precisely those of
(\ref{2.4}). \qed

\section{Proof of Theorem 1}
\setcounter{equation}{0}

We need to quote some results from [DM2]. First, following [Z] and
[DM2], we say that the VOA $V$ satisfies the {\em Virasoro condition}
if $V$ is a sum of highest weight modules for the Virasoro algebra
generated by the components of
$Y(\o,z);$ we say that $V$ satisfies the $C_2$ {\em condition} if
$V/C_2(V)$ is finite-dimensional, where $C_2(V)$ is spanned by
$u_{-2}v$ for $u,v\in V.$ What is important for us is that because
$\V$ contains $L$ as a rational subalgebra, $\V$ satisfies
the $C_2$ condition as well as the Virasoro condition. See [Z] for more
information on this point.

\begin{thm}\label{t3.1} [DM2] Suppose that $V$ is a holomorphic VOA satisfying
both
Virasoro and $C_2$ conditions. Let $g\in \Aut(V)$ have finite order. Then the
following hold:

(i) $V$ has at most one simple $g$-twisted module.

(ii) $V$ has at least one weak simple $g$-twisted module.
\end{thm}

So to prove parts (i) and (ii) of
Theorem \ref{t1}, it is enough to prove that a weak simple
$g$-twisted $\V$-module for $g$ of type $2A, 2B$ or $4A$ has finite-dimensional
homogeneous spaces. Let $M$ be such a module.

Let $W=(\V)^{\<g\>}.$ Then $M$ is an ordinary weak $W$-module, in fact it is
the sum of a finite number of simple weak $W$-modules (see [DM3], for example).
So it suffices to prove that a weak simple $W$-module $N,$ say,
has finite-dimensional homogeneous spaces.

It is shown in [DMZ]
that $L(\frac{1}{2},0)$ is a rational VOA with three irreducible modules
$L(\frac{1}{2},h_i),$ $h_i=0,1/2,1/16$ which are exactly the highest weight
unitary representations of $Vir$ with central charge $1/2.$
Moreover all weak modules are ordinary modules. The characters of these modules
are as follows: as well as $L(\frac{1}{2},0)$ given in (\ref{2.3}) we have
\begin{equation}\label{3.1}
\ba{c}
\displaystyle{{\rm
ch}_qL(\frac{1}{2},\frac{1}{2})=\frac{1}{2}\left(\prod_{n=1}^{\infty}(1+q^{n-1/2})-
\prod_{n=1}^{\infty}(1-q^{n-1/2})\right)} \\[10pt]
\displaystyle{{\rm
ch}_qL(\frac{1}{2},\frac{1}{16})=q^{1/16}\prod_{n=1}^{\infty}(1+q^{n}).}
\ea
\end{equation}
The fusion rules are as follows:
$$L(\frac{1}{2},\frac{1}{2})\times L(\frac{1}{2},\frac{1}{2})
=L(\frac{1}{2},0)$$
\begin{equation}
L(\frac{1}{2},\frac{1}{2})\times L(\frac{1}{2},\frac{1}{16})=
L(\frac{1}{2},\frac{1}{16})\label{3.2}
\end{equation}
$$L(\frac{1}{2},\frac{1}{16})\times L(\frac{1}{2},\frac{1}{16})
=L(\frac{1}{2},0)+L(\frac{1}{2},\frac{1}{2})$$
where the ``product'' $\times$ is commutative and
$L(\frac{1}{2},0)$ is the identity.

{}From these facts we conclude (cf. [FHL]) that $L$ has just $3^{48}$
simple modules, namely
\begin{equation}\label{3.3}
L(h_1,...,h_{48})=L(\frac{1}{2},h_1)\otimes\cdots\otimes L(\frac{1}{2},h_{48})
\end{equation}
where $h_i\in\{0,1/2,1/16\}.$ The corresponding fusion rules (see Proposition
2.10 of [DMZ]) are
\begin{equation}\label{3.3'}
L(h_1,...,h_{48})\times L(h_1',...,h_{48}')
=(L(\frac{1}{2},h_1)\times L(\frac{1}{2},h_1'))\otimes\cdots\otimes
(L(\frac{1}{2},h_{48})\times L(\frac{1}{2},h_{48}'))
\end{equation}

The fusion rules can be interpreted in terms of a tensor product $\boxtimes$
of modules for vertex operator algebras (see [HL] and [L]).
Then the product $\times$ in (\ref{3.2}) and
(\ref{3.3'}) can be replaced by the tensor product $\boxtimes$ to get the
tensor product decomposition. We refer the reader to [L] for more details.

Now if $0\ne n\in N$ is a highest weight vector
with highest weight $(h_1,...,h_{48})$
for the Virasoro
algebra generated by the components of
$Y(\o_i,z),$ $i=1,...,48,$ then the $L$-module generated by
$n$ is isomorphic to $L(h_1,...,h_{48}).$ Let $N_{h_1',...,h_{48}'}$ be
the multiplicity of $L(h_1',...,h_{48}')$ in $W.$ Then as an $L$-module $W$
has decomposition
$$W=\oplus_{h_i'\in\{0,1/2,1/16\}}N_{h_1',...,h_{48}'}L(h_1',...,h_{48}').$$
Then we have the tensor product of $L$-modules
$$W\boxtimes L(h_1,...,h_{48})=\oplus_{h_i'}N_{h_1',...,h_{48}'}
L(h_1',...,h'_{48})\boxtimes L(h_1,...,h_{48})$$
which is an ordinary $L$-module with finite-dimensional homogeneous subspaces.
{}From Proposition 4.1 of [DM3] we see that $N$ is spanned by
$w_mn$ for $w\in W$ and $m\in\Z.$ Now using the universal property of tensor
product, we conclude that $N$ is a submodule of $W\boxtimes L(h_1,...,h_{48})$
as $L$-modules. Thus each homogeneous subspace of $N$ is finite-dimensional.
This completes the proof of parts (i) and (ii) of Theorem 1.

If $H=2Baby$ is the centralizer of $g=2A$ in $\M$ then, as explained
in [DM1] and [DM4], the {\em uniqueness} of $\V(2A)$ yields a
projective representation of $H$ on $\V(2A).$ This must be an ordinary
and faithful
representation since $H^2(Baby,\C^*)\simeq \Z_2.$ Thus the group of extended
automorphisms of $\V(2A)$ is precisely $\C^*\times H,$ as claimed in the
introduction.

\section{Top level of $\V(2A)$}
\setcounter{equation}{0}

First we review some further results from [DM2]. Under the assumptions of the
Virasoro condition and the
$C_2$ condition, which we know hold for $\V,$ together with
the complete reducibility of $\V$-modules [D],
 the modular-invariance
properties established in [DM2] may be stated as follows:
\begin{equation}\label{4.1}
Z(g,h,\gamma\tau)=\sigma(\gamma^{-1},g,h)Z((g,h)\gamma,\tau).
\end{equation}
Here $(g,h)$ is a pair of elements which generate a cyclic group,
$Z(g,h,\tau)$ is a function in the so-called $(g,h)$-conformal block,
$\gamma\in SL(2,\Z)$ and $\sigma(\gamma^{-1},g,h)$ is a nonzero
constant. The notation $(g,h)\gamma$ is the action of $SL(2,\Z)$ on
pairs of commuting elements:
\begin{equation}\label{4.2}
(g,h)\left(
\ba{cc}
a & b\\
c & d
\ea
\right)=(g^ah^c,g^bh^d).
\end{equation}
What is important is that if $g$ is either $1A$ or one of $2A,2B$ or
$4A,$ so that Theorem \ref{t1} applies, then $qZ(g,h,\tau)$ is
precisely the graded trace of $h$ on the $g$-twisted sector
$\V(g).$ In particular, taking $g=1,$
$h=2A$ and $\gamma=\left(
\ba{cc}
0 & -1\\
1 & 0
\ea
\right)=S$ in (\ref{4.1}) yields
\begin{equation}\label{4.3}
Z(1,2A,S\tau)=\sigma(S^{-1},1,2A)Z(2A,1,\tau).
\end{equation}

By Borcherds' theorem [B2], $Z(1,2A,\tau)$ is a hauptmodul, specifically the
one denoted by $2+$ in [CN]. This means that $Z(1,2A,\tau)$ is invariant
under the Fricke involution $W_2=\left(
\ba{cc}
0 & -1\\
2 & 0
\ea
\right)=S\left(
\ba{cc}
2 & 0\\
0 & 1
\ea
\right),$
so that (\ref{4.3}) yields
$$Z(1,2A,\tau/2)=Z(1,2A,W_2(\tau/2))=Z(1,2A,S\tau)=\sigma(S^{-1},1,2A)Z(2A,1,\tau).$$
Now $Z(1,2A,\tau/2)=q^{-1/2}+O(q^{1/2}),$  and  $qZ(2A,1,\tau)$ is the
$q$-character of $\V(2A).$ We conclude that
\begin{lem}\label{l4.1}
The $q$-character of $\V(2A)$ has the form
$$k(q^{1/2}+O(q^{3/2}))$$
for some constant $k\ne 0.$
\end{lem}

We will establish
\begin{thm}\label{t4.2}
The constant $k$ is equal to 1.
\end{thm}

In the following we let $M=\V(2A)$ with $M_{1/2}$ the top level of $M.$
Now $M_{1/2}$ is spanned by highest weight vectors corresponding to
simple $L$-modules $L(h_1',...,h_{48}')\subset M$ satisfying
$\sum_{i}h_i'=1/2.$  Then $M^{1}=\sum_{n\in \Z}M_{n+\frac{1}{2}}$
and $M^{0}=\sum_{n\in \Z}M_n$ are irreducible $W$-modules where
$W=(\V)^{\<2A\>}$ (see Theorem 6.1 of [DM4]). Since we are mainly concerned
with the top level $M_{1/2}$ we will only pay attention to the
$W$-module $M^1.$ We have already explained that $W$ is the sum
of simple $L$-modules $W=W_1+\cdots + W_t$ and $M^1$ is spanned by
$w_lm$ for $w\in W,$ $l\in \Z,$ and fixed $0\ne m\in M_{1/2}.$ So if
we choose $m$ to be a highest weight vector in the simple
$L$-module $N=L(h_1',...,h_{48}')\subset M^1$ we get information about
$M^1$ by considering the fusion rules $W_i\times N$ or the tensor product
$W_i\boxtimes N.$

In particular, we can bound $\dim M_{1/2}$ by estimating how many
indices $i$ are such that $W_i\boxtimes N$ contains a simple $L$-module
of highest weight $1/2.$ Let
$W_i=L(h_1,...,h_{48})$ and assume that $W_i\boxtimes N\supset
L(k_1,...,k_{48})$ with $\sum_ik_i=1/2.$

\begin{lem}\label{l4.3} If $N=L(\frac{1}{2},0^{47})$ with `1/2' in any position
then
Theorem \ref{t4.2} holds.
\end{lem}

\pf Without loss take `1/2' in the first position. Then we have
$k_j=h_j,$ $2\leq j\leq 48,$ and $k_1=1/2$ if $h_1=0;$
$k_1=0$ if $h_1=1/2;$ $k_1=1/16$ if $h_1=1/16.$ This follows from
the fusion rules (\ref{3.2}).

Since $\sum_jh_j$ is a nonnegative integer not equal to 1 (since $\V_n=0$
for $n<0$ or $n=1$), the condition $\sum_{i}k_i=1/2$ forces $h_1=\cdots
h_{48}=0,$ that is $W_i=L.$ Since $L$ has multiplicity 1 in $\V$ and
$L\boxtimes N=N,$
the lemma follows. \ \ $\Box$
\bigskip

{}From now on  we assume that all highest
weight vectors in $M_{1/2}$ generate $L$-modules
of type $L((\frac{1}{16})^{8},0^{40})$ for some distribution of the
$1/16$'s. For now we take $N=L((\frac{1}{16})^{8}|0^{40}),$ meaning that
the $1/16$'s are in the first 8 coordinate positions (this is purely
for notational convenience).
\begin{lem}\label{l4.4}
If $W_i$ is such that $W_i\times N\supset N'=L(k_1,...,k_{48})$ with
$\sum_ik_i=1/2$ then one of the following holds:

(a) $W_i=L(0,(\frac{1}{2})^{3},(\frac{1}{16})^{4} | (\frac{1}{16})^{4},
0^{36})$ and $N'=L((\frac{1}{16})^{4},0^4|(\frac{1}{16})^{4},0^{36});$

(b) $W_i=L((\frac{1}{2})^{8}|0^{40})$ and $N'=N;$

(c) $W_i=L((\frac{1}{2})^{6},0^2|0^{40})$ and $N'=N;$

(d) $W_i=L((\frac{1}{2})^{4},0^4|0^{40})$ and $N'=N;$

(e) $W_i=L$ and $N'=N.$

\noindent(No specific ordering of the first 8 coordinates or the last 40
entries is implied.)
\end{lem}

\pf Let $W_i=(0^a,(\frac{1}{2})^{b},(\frac{1}{16})^{c}|h_9,...,h_{48}).$ Using
the fusion rules (\ref{3.2}), we see that
$N'=((\frac{1}{16})^{a+b},0^{c-d},(\frac{1}{2})^{d}|h_9,...,h_{48})$ for some
$0\leq d\leq c.$

We have, setting $s=\sum_{i=9}^{48}h_i,$ that
$$\frac{a+b}{16}+\frac{d}{2}+s=\frac{1}{2}$$
$$\frac{b}{2}+\frac{c}{16}+s=0\ \ {\rm or}\ \geq 2\ {\rm, lies\ in} \ \Z$$
$$a+b+c=8.$$
It follows that $c\equiv 0$ ($\mod\,4$) and $d=0.$ If $c=8$ then $a+b=0,$
so $s=1/2$ and $c/16+s=1,$ contradiction. So $c=0$ or 4. If $c=4$ then
$a+b=4,$ $s=1/4,$ and then $b=3.$ If $c=0$ then $a+b=8, s=0,
b=4,6$ or 8. The lemma follows. \ \ $\Box$

\begin{lem}\label{ladd} In the notation of Lemma \ref{l4.4}, if
$W_i \times N=N'\subset M^1$ then $w_{\wt\,w-1}m\ne 0$ where $w$ is
a nonzero highest weight vector of $W_i.$ Moreover $v_{\wt v-1}m$ is a
scalar multiple of $w_{\wt w-1}m$ for any homogeneous $v\in W_i.$
\end{lem}

\pf We need to use results on the Zhu algebra $A(V)$ and its bimodule
$A(W_i)$ (which is a quotient of $W_i$ by a subspace)
to prove this result. We refer the reader to [Z] and [FZ] for
the definitions. Let $m'$ be a nonzero highest weight vector of $N'.$ By
Theorem 1.5.3 of [FZ] and Theorem 4.2.4 of [L], $A(W^i)\otimes_{A(L)}\C m$
is isomorphic to $\C m'$ as $A(L)$-modules under the map
$\bar v\otimes m\to v_{\wt v-1}m$ where $\bar v$ is the image of $v\in W_i$
in $A(W_i).$ By Lemma 2.8 and Proposition 3.3 of [DMZ], $A(L)$ is isomorphic to
the associative commutative algebra $\C[t_j|j=1,...,48]/I$ where $I$
is the ideal generated by $t_j(t_j-\frac{1}{2})(t_j-\frac{1}{16}).$
By Lemma 2.9, Propositions 3.1 and 3.4 of [DMZ], $A(W_i)$ is isomorphic to
$\C[x_j,y_j|j=1,...,48]/I_i$ where $I_i$ is a certain ideal of
$\C[x_j,y_j|j=1,...,48]$ and the left and right actions of $t_j+I$ on $A(W_i)$
are multiplications by $x_j$ and $y_j$ respectively. Moreover, under this
identification, $\bar w$ is mapped to $1+I_i.$ Since $A(W_i)$ is generated
by $\bar w$ as $A(L)$-bimodule, the lemma follows immediately.
\qed

\begin{lem}\label{l4.5} $M_{1/2}$ is a trivial module for $H=2Baby.$
\end{lem}

\pf The minimal nontrivial representation of $H$ has degree 4371 [Cal] so
it suffices to show that $\dim M_{1/2}$ is less than this. Now the
multiplicity of $L(h_1,...,h_{48})$ in $\V$ is less than or equal to 1
if all $h_i\in\{0,1/2\}$ by Proposition 5.1 of [DMZ]. So by Lemma
\ref{l4.4} we conclude that the multiplicity of
$N=L((\frac{1}{16})^8|0^{40})$ in $M^1$ is at most $ {8\choose 0}
+{8\choose 4}+{8\choose 6}+{8\choose 8}=100.$ All other simple
$L$-submodules of $M^1$ of weight $1/2$ arise from the action of
$L(0,(\frac{1}{2})^{3},(\frac{1}{16})^{4} | (\frac{1}{16})^{4},
0^{36})\subset \V_2,$ and the multiplicity\footnote{Some of the multiplicities
stated in Theorem 6.5 of [DMZ] are inaccurate. See [Ho] for corrections.}
$\mu$ of {\em all} modules of
this type in $\V$ is given in Theorem 6.5 [DMZ]. Since $\mu=24\cdot 2^6,$
the lemma follows.\ \ $\Box$

\begin{lem}\label{l4.6} $M_{1/2}$ is contained in the span
of $v_{\wt v-1}m$ for $v\in \V_n$ with $n\leq 4.$
\end{lem}

\pf First we observe from Lemma \ref{l4.4} that the highest weight vectors
in $W_i$ have weights 0,2,3 or 4. By Corollary 4.2 of [DM3], $M_{1/2}$ is
spanned by $v_{\wt\,v-1}m$ for $v\in W_i$  with $W_i$ occurring
in Lemma \ref{l4.4}.  So
it is enough to show that the span of $v_{\wt\,v-1}m$ for $v\in W_i$
is contained in $w_{\wt\,w-1}m$ where $w$ is a nonzero highest weight vector
of $W_i.$ This follows from Lemma \ref{ladd}. \qed
\bigskip

Now we can complete the proof of Theorem \ref{t4.2}. We may represent
the conclusion of Lemma \ref{l4.6}
by the containment  $(\oplus_{n=0}^4\V_n)m\supset M_{1/2}.$
Since $M_{1/2}$ is a trivial $H$-module by Lemma \ref{l4.5}, it follows
that $M_{1/2}\subset (\oplus_{n=0}^4(\V_n)^H)m.$

Now use Lemma \ref{l2.5} to see that $M_{1/2}$ lies in the space
spanned by $u_{\wt\,u-1}m$ for $u\in L(\frac{1}{2},0)\otimes L(\frac{47}{2},0)$
homogeneous. But $m$ is an eigenvector for such operators. So we get
$M_{1/2}=\C m,$ as required. \qed
\bigskip

We are now ready for the proof of Theorem \ref{t1} (iii). Let $M$ be an
irreducible $g$-twisted $\V$-module so that the top level of $M$
is $\C m.$ Let $N$ be the $L$-submodule of $M$ generated by $m.$ Then
by Lemmas \ref{l4.3} and \ref{l4.4}, $N$ is either $L(\frac{1}{2},0^{47})$
or $L((\frac{1}{16})^{8}|0^{40}).$ Then $W_i\times N=N$ if, and only, if
in the former case $W_i=L,$ and in the latter case $W_i$ appears in
(b)-(e) of Lemma \ref{l4.4}. In either case, since $M$ is simple
we see that the subspace of $N$ spanned
by $v_{l}m$ for $v\in W_i$ and $l\in\Z$ is $N$ (by Proposition 11.9
of [DL]). Again by Lemma \ref{ladd}, $w_{\wt\,w-1}m\ne 0.$

Let $X$ be a $g$-twisted $\V$-module. Then $X$ has a finite composition
series (see [DLM]). Using induction on the number of composition-factors
of $X,$ we
only need to prove that $X$ is completely reducible if $X$ has two
factors. It is shown in [DLM] that
 $X$ is a completely reducible $g$-twisted $\V$-module
if, and only, if $X_{\frac{1}{2}}$ is a semisimple $A_g(\V)$-module
via the action $v_{\wt\,v-1}$ for homogeneous $v\in\V.$ If
$N=L(\frac{1}{2},0^{47})$ it is clear that $X_{\frac{1}{2}}$ is
a semisimple $A_g(\V)$-module from the discussion above. Here $A_g(\V)$
is the twisted Zhu algebra as defined and used in [DLM].

Now we assume that $N=L((\frac{1}{16})^{8}|0^{40}).$ Let
$X_{\frac{1}{2}}=\C x_1+\C x_2$ such that $x_1$ generates an
irreducible $g$-twisted $\V$-module $X^1$ (which is necessarily
isomorphic to $M$). Using the associativity of vertex operators (see
the proof of Proposition 4.1 of [DM3], for example) we see that
$w_{\wt\,w-1}w_{\wt\,w-1}x_2$ is a nonzero multiple of $x_2.$ Thus
$w_{\wt\,w-1}$ acts semisimply on $X_{\frac{1}{2}},$ hence acts as a
scalar. Thus for all homogeneous $v\in W_i$ the action of
$v_{\wt\,v-1}$ on $X_{\frac{1}{2}}$ is semisimple.  Note that the
image of $A_g(\V)$ in $\End(X_{\frac{1}{2}})$ is a subalgebra of
dimension less than or equal to 2. So we conclude that $X_{\frac{1}{2}}$
is indeed a semisimple $A_g(\V)$-module.

\section{Proof of Theorem 2}
\setcounter{equation}{0}
Let $g\in\M$ be of type $2A$ and let $h$ be an element of $C_{\M}(g)$ of
{\em odd} order $N.$ We use the notation of (\ref{4.1})-(\ref{4.2}).
\begin{lem}\label{l5.1}
Let $F\leq SL(2,\R)$ be the fixing group of $Z(1,gh,\tau).$ Then $F$ contains
the Atkin-Lehner involution $W_2.$
\end{lem}

\pf By Borcherds' theorem [B2], each $Z(1,x,\tau)$ for $x\in\M$ is a hauptmodul
on a discrete group $F=F_x\leq SL(2,\R).$ Then $F$ is precisely the group
conjectured in [CN], Table 2. This has been established by the work of several
authors; see [CN] and [F] for further references.

It is a fact that if $x=gh$ then $W_2$ always lies in $F.$ In fact, the group
$F$ is of the form $2N+,$ or $2N+2$ in all but one case. The exception is the
element $30F,$ where the group is 30+2, 15,30. In any case, the lemma follows.
\ \ $\Box$

\bigskip
Now we may take   $ W_2=\left(
\ba{cc}
a & b\\
cN & 2d
\ea
\right)\left(
\ba{cc}
2 & 0\\
0 & 1
\ea
\right)$ where $\gamma=\left(
\ba{cc}
a & b\\
cN & 2d
\ea
\right)\in SL(2,\Z).$ Then we have from (\ref{4.1}) that
\begin{eqnarray*}
& &Z(1,gh,\tau/2)=Z(1,gh,W_2(\tau/2))\\
& &\ \ =\sigma Z((1,gh)\gamma,\tau)\\
& &\ \ =\sigma Z((gh)^{cN},(gh)^{2d},\tau)\\
& &\ \ =\sigma Z(g,h',\tau).
\end{eqnarray*}
Here, $\sigma$ is a constant and $h'=h^{a'}$ where $aa'\equiv 1$ ($\mod\,N$).
But $Z(1,gh,\tau/2)=q^{-1/2}+\cdots $ and $Z(g,h',\tau)=q^{-1/2}+\cdots $
by Theorem \ref{t4.2}. So $\sigma=1.$ Now Theorem \ref{t2}, part (i) follows
immediately.

As for part (ii) of Theorem \ref{t2}, if we now take $h\in\M$ such that
$g\in\<h\>$ then we can find $\gamma\in SL(2,\Z)$ such that
$(g,h)=(1,h)\gamma.$ Then (\ref{4.1}) yields
\begin{equation}\label{5.1}
Z(g,h,\tau)=\sigma Z(1,h,\gamma\tau)
\end{equation}
for some constant $\sigma.$ As $Z(1,h,\tau)$ is a hauptmodul by Borcherds'
theorem  then so is $Z(g,h,\tau)$ by (\ref{5.1}).

Now let us consider $Z(g,h,\tau)$ where $g,h$ and $gh$ are all of type
$2A.$ From [DM2] we know that $Z(g,h,\tau)$ is holomorphic on the upper
half-plane and meromorphic at the cusps. Furthermore (\ref{4.1}) still
holds because of Theorem \ref{t1} (iii). We see that if $\gamma\in SL(2,\Z)$
then
$$Z(g,h,\gamma\tau)=\sigma(\gamma)Z(g,h,\tau)$$
for some constant $\sigma(\gamma).$ So $\sigma$ is a character of
$SL(2,\Z).$ As $T=\left(
\ba{cc}
1 & 1\\
0 & 1
\ea\right)$ covers the abelianization of $SL(2,\Z)$ then the kernel
$K$ of $\sigma$ is the subgroup of $SL(2,\Z)$ of index 2. As $i\infty$
is the unique cusp for $K$ we see easily that $K$ is indeed of genus zero
with hauptmodul $Z(g,h,\tau).$ It is in fact
the function denoted
$t_{2/2}=\sqrt{J(q)-984}=q^{-1/2}-492q^{1/2}-22590q^{3/2}+\cdots$ in [CN]. To
see this, note that from
Theorem \ref{t2} (i) with $h=1,$ combined
with tables in [CN] and the results of [B2], we see that
$Z(g,1,\tau)=q^{-1/2}+4372q^{1/2}+\cdots.$ This tells us that the
weight $\frac{3}{2}$ subspace of $\V(2A)$ is the module $1\oplus 4371$
for $2Baby$ (see [Cal]), on which $h$ has trace $-492$ (ibid). Thus
$Z(g,h,\tau)=q^{-1/2}-492q^{1/2}+\cdots $ as claimed.

\section{Final comments}
\setcounter{equation}{0}

We have rather ignored the twisted sectors $\V(2B)$ and $\V(4A).$ As
we have said, there is a construction of $\V(2B)$ in [Hu], and
its existence also follows from [DM2] without the necessity of the effort
we needed to understand $\V(2A).$ Huang also constructs the
$2B$-{\em orbifold}, i.e., puts an abelian intertwining algebra
structure [DL] on $\V\oplus \V(2B).$ This is closely related to the
construction of $\V$ in [FLM].

Concerning the extended automorphism group it is known [G] that the
centralizer $C$ of $2B$ in $\M$ is a non-split extension of $\cdot 1$
(largest simple Conway group) by the extra-special group
$2^{1+24}_+.$ Furthermore (loc.cit.) $H^2(C,\C^*)\simeq \Z_2.$ In fact if
$\cdot 0$ is the 2-fold cover of $\cdot 1$ (=automorphism group of the
Leech lattice) there is a diagram
$$\begin{array}{ccccccccc}
1 &\to &2_+^{1+24}&\to &C&\to &\cdot 1&\to & 1\\
  & & \|& & \uparrow & & \uparrow & &\\
1 &\to &2_+^{1+24}&\to &\hat C&\to &\cdot 0&\to & 1
\end{array}
$$
and $\hat C$ is the universal central extension of $C.$

Griess also shows that $\hat C$ has a simple module of degree
$2^{12},$ whereas $C$ has no such representation. The smallest
faithful irreducible representation for $C$ has dimension $24\cdot
2^{12}.$

Now the weight $3/2$ subspace of $\V(2B)$ has dimension $2^{12}$ (see [Hu] and
[DM2]) and there
is a projective representation of $C$ on $\V(2B)$ by [DM1] and [DM4]. The
conclusion is thus
\begin{lem}\label{l6.1}
The extended automorphism group of $\V(2B)$ is a non-split extension
$$1\to \C^*\to \Aut^e(\V(2B))\to C\to 1.$$
It has commutator subgroup isomorphic to $\hat C.$
\end{lem}

\begin{thm}\label{t6.2}
Let $h\in \M$ be such that $\<h\>$ contains $g=2B$ or $4A.$ Then
$Z(g,h,\tau)$ is a hauptmodul.
\end{thm}

\pf Same as the proof of Theorem \ref{t2} (ii).

\vspace{0.5 cm}

Finally, we consider the $q$-characters of $\V(2B)$ and $\V(4A).$
Note by [CN] or [FLM] that we have for $g\in \M$ of type $2B$ that
$$Z(1,g,\tau)=24+q^{-1}\prod_{n\,odd}(1-q^n)^{24}
=24+\frac{\eta(\tau)^{24}}{\eta(2\tau)^{24}}$$
where $\eta(\tau)$ is the Dedekind eta-function. By (\ref{4.1}) and
the transformation law for $\eta(\tau)$ [K] we get the $q$-character
of $\V(2B)$ equal to
$$Z(g,1,\tau)=\sigma
Z(1,g,S\tau)=\sigma\left\{24+\frac{2^{12}\eta(\tau)^{24}}{\eta(\tau/2)^{24}}\right\}.$$
In fact $\sigma=1,$ as we know from [Hu] and [DM2].

Similarly for $t$ of type $4A$ we have
$$Z(1,t,\tau)=-24+\frac{\eta(2\tau)^{48}}{\eta(\tau)^{24}\eta(4\tau)^{24}}$$
which leads to
$$Z(t,1,\tau)=\sigma\left\{-24+\frac{\eta(\tau/2)^{48}}{\eta(\tau)^{24}\eta(\tau/4)^{24}}\right\}.$$
Presumably $\sigma=1$ in this case too, but a proof would be more
difficult than that for $2A$ given above.

\end{document}